\documentclass[11pt]{extarticle}
\usepackage[cp1251]{inputenc}
\usepackage[russian]{babel}
\usepackage{amsmath}
\usepackage{amssymb}
\usepackage[dvips]{graphicx}
\usepackage{wrapfig}
\usepackage{floatflt}
\usepackage{cmap}

\newtheorem{theorem}{Lemma}

\newcommand{\Name}[1]{\section{\large #1}}

\begin{document}
\renewcommand{\abstractname}{Abstract}
\renewcommand{\refname}{References}
\title{One-dimensional minimal fillings \\ with negative edge weights}
\author{A.\,O.~Ivanov\and Z.\,N.~Ovsyannikov \and N.\,P.~Strelkova\and A.\,A.~Tuzhilin}
\date{}

\maketitle

\begin{abstract}

Ivanov and Tuzhilin started an investigation of a particular case
of Gromov Minimal Fillings problem (generalized to the case of
stratified manifolds). Weighted graphs with non-negative weight
function were used as minimal fillings of finite metric spaces. In
the present paper we introduce generalized minimal fillings, i.e.
minimal fillings where the weight function is not necessarily
non-negative. We prove that for any finite metric space its
minimal filling has the minimum weight in the class of all
generalized fillings of the space.

\smallskip

\end{abstract}

\Name{ Introduction}

\noindent The problem concerning minimal fillings of finite metric
spaces was posed by Ivanov and Tuzhilin in~\cite{it}. This problem
arose in connection with two classical problems, the Steiner tree
problem
 and the problem on minimal fillings of
smooth manifolds (posed by M.~Gromov, see~\cite{gromov}).

The objective is to find a weighted graph of minimum weight among
all weighted graphs joining the points of a given finite metric
space provided that for any two points in the metric space the
distance between them is not greater than the weight of the shortest
path connecting them in the graph.  Now we give the precise
definitions from~\cite{it}.

Let $\mathcal{M} = (M,\rho)$ be a finite pseudometric space (a
pseudometric space is a generalized metric space in which the
distance between two distinct points can be zero), let  $G = (V,E)$
be a connected graph joining $M$ (by definition a connected graph
$G$ {\it joins} $M$ if $M\subset V$), let $w\colon
E\rightarrow\mathbb{R}_+$ be a function from the graph edge set to
the set of non-negative real numbers, which is called the {\it
weight function} of  the {\it  weighted graph } $\mathcal{G}=(G,w)$.
The sum $\sum_{e\in E} w(e)$ of all edge weights of $\cal G$ is
called {\it the weight of the graph} and is denoted by $w({\cal
G})$. We define a distance function $d_w$ on the set  $V$ by saying
that the distance between two points is equal to  the weight of a
shortest path connecting these two points in the graph.

The weighted graph  ${\cal G}=(G,w)$ is called a {\it filling} of
the
 space $\mathcal{M}$ if for any two points $p,q\in M$ the
 condition
$\rho(p,q)\leqslant d_w(p,q)$ holds. In this case the graph $G$ is
called the {\it type} of this  filling. The number
$\mathrm{mf}(\mathcal{M})=\inf w(\mathcal{G})$, where infimum is
taken over all fillings $\mathcal{G}$ of the space $\mathcal{M}$ is
called the {\it minimal filling weight}, while a filling
$\mathcal{G}$ provided $w(\mathcal{G})=\mathrm{mf}(\mathcal{M})$ is
called a {\it minimal filling} of the space $\mathcal{M}$.

In paper~\cite{it}, only weighted graphs with non-negative weight
function were used as fillings. In the present paper, we don't
forbid negative weight edges and are looking for the minimum graph
weight among fillings with arbitrary, not necessarily non-negative
weight function. Such fillings will be called {\it generalized}.

The generalized minimal fillings have most of the properties proved
in~\cite{it} for non-generalized minimal fillings. Some properties
will be proved in section~3.

The main result of the present paper is the surprising fact that for
any metric space its minimal filling weight (i.e.  the minimum of
weight on the set of fillings with non-negative weight function)
equals the minimum of weight on the larger set of generalized
fillings (i.e. fillings with arbitrary weight function).

This result simplifies the original problem of finding the minimal
filling weight since now one does not have to verify whether the
edge weights are non-negative. Using this result, A.\,Yu.~Eremin
in~\cite{eremin} has proved a formula for calculating minimal
filling weight of a finite metric space.

The authors are grateful to academician A.\,T.~Fomenko for
constant attention to this work, and are also grateful to  the
participants of the  seminar "Minimal networks"\ in Mechanics and
Mathematics Department at MSU  for useful discussions.

\Name{ Definitions}

\noindent We have already given the definitions of a filling and a
minimal filling of a finite pseudometric space in the Introduction.

Let a graph $G=(V,E)$  and an embedding $M\subset V$ be fixed.
Consider $\inf w(\mathcal{G})$ taken over all the fillings $\cal G$
of $\cal M$  having the fixed type~$G$. This infimum is called the
{\it minimal parametric filling weight of the fixed type~$G$} and is
denoted by $\mathrm{mpf}(\mathcal{M},G)$, while any filling of
type~$G$ on which the infimum is reached  is called a {\it minimal
parametric filling of type~$G$}, see~\cite{it}.

So-defined minimal and minimal parametric fillings exist for any
finite pseudometric space (\cite[theorem 2.1]{it}).

In a graph $G$ joining the given finite set $M$, the vertices that
correspond to the points from $M$ will be called {\it boundary}, and
the remaining vertices in the graph will be called {\it interior}.
The graph edge connecting  vertices $u$ and $v$ will be denoted
by~$uv$.

Now we modify the definitions of fillings given above by allowing
the edge weights be negative.

 Let $G$ be a graph and $w\colon
E\rightarrow \mathbb{R}$ an arbitrary function. Then the pair
$(G,w)= (V,E,w)$ is called a {\it generalized weighted graph}. We
define $d_w\colon V\times V \rightarrow \mathbb{R}$ by saying that
$d_w(u,v)$ is the smallest among weights of simple (i.e.,
non-self-intersecting) paths connecting $u$ and $v$. The function
$d_w$ is not necessarily non-negative and may violate the triangle
inequality.

A generalized weighted graph $\cal G$ joining $M$ is called a {\it
generalized filling} of a finite pseudometric space ${\cal M}
=(M,\rho)$ if for any $u,v\in M$ holds: $\rho(u,v)\leqslant
d_w(u,v)$.

Fix a graph $G=(V,E)$  and an embedding $M\subset V$. The value
$\inf w({\cal G})$, where the infimum is taken over all generalized
fillings ${\cal G}$ of $\cal M$ of the fixed type $G$ is called the
{\it generalized parametric minimal filling weight of type $G$} and
is denoted by $\mathrm{mpf}_{-}(\mathcal{M},G)$. A generalized
filling $\cal G$ of type $G$ is called a {\it generalized parametric
minimal filling of type $G$} if $w({\cal
G})=\mathrm{mpf}_{-}(\mathcal{M},G)$.

The value $\inf {\rm mpf}_{-}({\cal M}, G)$, where the infimum is
taken over all trees $G$  joining $M$ and satisfying the condition
that any vertex $v\in G$ of degree $1$ is boundary (i.e., $v\in M$)
is called the {\it generalized minimal filling weight of type $G$}
and is denoted by $\mathrm{mf}_{-}(\mathcal{M})$.  A generalized
filling $\cal G$ is called a {\it generalized minimal filling of
type $G$} if $w({\cal G})=\mathrm{mf}_{-}(\mathcal{M})$.

\medskip

\noindent{\bf Remark 1.} There is a similar definition for
non-generalized minimal fillings: ${\rm mf}({\cal M}) = \inf{\rm
mpf}({\cal M},G)$ where there is no difference whether to take the
infimum over trees or over arbitrary graphs  $G$ joining  $M$. The
definition in both cases is equivalent to the definition of minimal
filling weight given above, see~\cite{it}. But the situation is
different in the case of generalized minimal fillings --- see
example~1.

\medskip

\noindent{\bf Remark 2.} One certainly obtains an equivalent
definition for generalized minimal filling if one takes $\inf {\rm
mpf}_{-}({\cal M},G)$ only over trees without interior vertices of
degree~$2$.

\medskip

\noindent{\bf Remark 3.} Suppose the embedding $M\subset V$ is
organized in a such way that there is an interior vertex of degree
$1$ in $G$. We claim that then ${\rm mpf}_{-}({\cal
M},G)=-\infty$. This easily follows from the fact that there are
no restrictions on the weight of the edge incident to this vertex
and so the weight of this edge can be made any negative number
regardless of what other edge weights are.

\medskip

\begin{wrapfigure}[12]{r}{60pt}
\includegraphics[scale = 0.6]{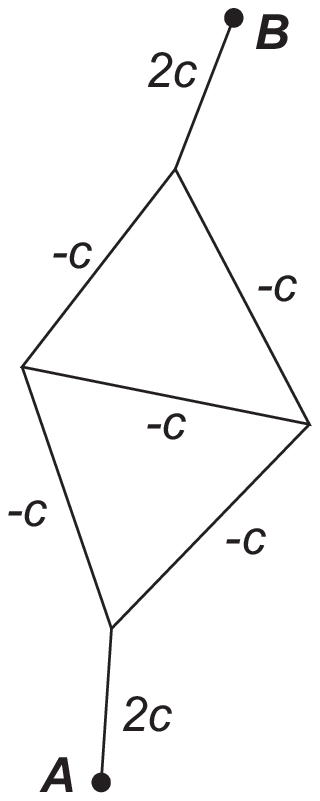}
\end{wrapfigure}

\noindent{\bf Example 1.} We give an example where ${\rm
mpf}_{-}({\cal M},G)=-\infty$ though for any edge in the graph there
exists a path with boundary endpoints that goes through this edge
(in contradiction to the cases like one considered above in which
there is an interior vertex of degree 1). Consider a metric space
consisting of two points  $M = \{A,B\}$ with distance $\rho(A,B)=1$.
The figure shows the graph $G$ joining $M$ and the edge weights $w$.
We assume $c>0$.

Now we have $d_w(A,A)=0, d_w(A,B)=c, d_w(B,B)=0$. Therefore for
any $c\geqslant1$ the graph $(G,w)$ is a generalized filling for
space~$\cal M$. But  $w(G)= -c$. A similar example can be easily
constructed for any finite metric space.

\smallskip

\noindent{\bf Example 2.} We give an example of a space~$\cal M$
and a tree~$G$ with the generalized parametric minimal filling
weight strictly less than the parametric minimal filling weight:
${\rm mpf}({\cal M},G)>{\rm mpf_{-}}({\cal M},G)$.

Consider $M = \{a,b,c,d\}$, $\rho(a,b)=\rho(c,d)=4$,
$\rho(b,c)=\rho(a,d)= \rho(a,c)=\rho(b,d)=3$.  Define graph~$G=
(V,E)$ as follows: $V= \{a,b,c,d,u,v\}$, $E=\{au,bu,uv,cv,dv\}$. We
claim that in this case the minimal parametric filling weight is~$8$
(the reason is that $w(au)+w(ub)\geqslant\rho(a,b)= 4,
w(cv)+w(vd)\geqslant \rho(c,d)=4,w(uv)\geqslant 0$ so ${\rm
mpf}({\cal M},G)\geqslant8$ and $8$ is the weight of the filling
with the weight function $w(uv)=0$, $w(e)=2$ for other edges). On
the other hand, it can be easily checked that the generalized
parametric minimal filling weight is the weight of the filling with
weights $w(uv)=-1$, $w(e)=2$ for other edges. Thus we have an
example where ${\rm mpf}({\cal M},G)=8>7={\rm mpf_{-}}({\cal M},G)$.

Clearly, every generalized filling is a filling. So ${\rm
mpf_{-}}({\cal M},G)\leqslant{\rm mpf}({\cal M},G)$ for any $G$ and
also ${\rm mf_{-}}({\cal M})\leqslant{\rm mf}({\cal M})$.

The main result of this paper is the following theorem, which is
proved in section~\ref{sec: the-main-theorem}.

\medskip

\noindent {\bf Theorem.} {\it Let ${\cal M}$ be a finite
pseudometric space. Then ${\rm mf_{-}}({\cal M})={\rm mf}({\cal
M})$.}

\Name{Properties of generalized minimal fillings}

\noindent We assume throughout this section that the type of any
generalized filling under consideration is a tree and all its
vertices of degree 1 or 2 are boundary (i.e. the vertices are
in~$M$).

A  {\it cyclic order}  on a finite set~$M$ with $n$ elements is an
arbitrary ordering of its elements into a list or in other words a
bijection $\pi\colon\mathbb{Z}_n \rightarrow M$. Consider a tree $G$
joining~$M$. For a given cyclic order~$\pi$ on~$M$ consider the
paths in~$G$ connecting~$\pi(k)$
 and~$\pi(k+1)$ where $k= 0,\ldots, n-1$. The collection of these paths is called
 a {\it tour through $G$ corresponding to the cyclic order $\pi$}.
 A cyclic order $\pi$ on $M$ is called {\it
planar with respect to $G$} if any edge in $G$ is contained in
exactly two paths from the tour corresponding to the cyclic order.
Note that in~[1] the term "tour" \ was used as a synonym of the
expression "planar order".

The following simple assertion will be useful for us.

\begin{theorem}\label{lem: planar cyclic order}
Suppose $M$ includes all vertices of degree 1 from $G$. Then there
exists a cyclic order on~$M$ that is planar with respect to~$G$.
\end{theorem}

\noindent {\bf Remark 4.} To construct such an order, one can embed
the tree~$G$ into the plane in arbitrary way and make a walk around
this tree in the plane and  denote by $\pi(k-1)$ the vertex number
$k$ from $M$ that one meets along the walk. The fact is that any
planar order can be obtained by such procedure. See~\cite[section
7]{it} for more details and definitions. We shall need only the
definition of planar order and lemma~\ref{lem: planar cyclic order}.

\begin{theorem}\label{lem: total weight is positive}
For any pseudometric space $\cal M$ its generalized minimal filling
weight is positive (excluding the trivial case of distance function
being identically zero in $\cal M$
--- in that case the generalized minimal filling weight is zero).
\end{theorem}
{\bf Proof.} Let $\cal G$ be a generalized filling for space  $\cal
M$ of type $G$ where $G=(V,E)$ is a tree with no interior vertices
of degree 1. Consider an arbitrary order~$\pi$ on $M$ that is planar
with respect to $G$. The definition of planar order implies that the
sum of all edge weights contained in the paths from the tour
corresponding to the planar order equals the doubled weight of the
filling. But the definition of generalized filling asserts that the
weight of each path is greater or equal to the distance between its
endpoints in space $\cal M$. Therefore:
$$w({\cal G})=
\frac{1}{2}\sum_{k=0}^{n-1}d_w(\pi(k),\pi(k+1))\geqslant\frac{1}{2}\sum_{k=0}^{n-1}\rho(\pi(k),\pi(k+1))>0\mbox{
.}$$

This completes the proof of lemma~\ref{lem: total weight is
positive}.

\begin{theorem}\label{lem: mpf- exists}
For any finite pseudometric space ${\cal M}=(M,\rho)$  and any graph
$G=(V,E)$  joining~$M$ there exists a generalized minimal parametric
filling of type~$G$.
\end{theorem}
{\bf Proof.} A weight function $w$ is determined by the string of
its values assigned to the graph edges, i.e. by a vector from
$\mathbb{R}^{|E|}$. The restrictions that distinguish the weight
functions corresponding to generalized fillings~$(G,w)$ from other
weight functions are linear non-strict inequalities. Therefore the
subset $\Omega$ in $\mathbb{R}^{|E|}$ defined by these restrictions
is  an intersection of finitely many closed half-spaces. The linear
function $w({\cal G})=x_1+\ldots + x_{|E|}$ is bounded below on
$\Omega$ by the assertion of lemma~\ref{lem: total weight is
positive}. Therefore this function attains its minimum on $\Omega$.
This completes the proof.

\medskip

A filling will be said to be a {\it binary tree} if the type of the
filling is a tree with all vertices of degree  $1$ or $3$, moreover,
with the property that all its vertices of degree $1$ are boundary.

Let a subset $F\subset E$ of the edge set of a graph $G= (V,E)$ be
fixed. Denote by $G_i=(V_i,F_i), i=1,\ldots,m$ the connected
components of the graph $(V,F)$. Consider the new graph
$G_F=(V_F,E_F)$ where $V_F = \{V_1,\ldots,V_m\}$ and two vertices
$V_i$ and $V_j$ are adjacent if and only if there exists an  edge
$v_iv_j\in E\setminus F$ such that $v_i\in V_i$ and $v_j\in V_j$.
The graph $G_F$ is said to be the quotient graph of $G$ by the set
$F$.

The operation inverse to taking a quotient is called {\it
splitting}. The graph $G_F$ is said to be {\it obtained} from the
graph $G$ {\it by splitting the vertex} $V_i$ if $ G_i$ is the only
component consisting of more than one vertex.  See~\cite{it} for
related definitions and properties of such operations.

\begin{theorem}\label{lem: mf- exists}
For any finite pseudometric space ${\cal M}=(M,\rho)$
 there exits a minimal filling.
 Moreover, there  exists a binary tree minimal filling.
\end{theorem}
{\bf Proof.} While searching for the minimal filling weight one can
take the infimum $\inf {\rm mpf}_{-}({\cal M},G)$ only over trees
$G$ with all vertices of degree $1$ and $2$ lying in $M$
(see~remark~2). The set of trees with no more than $|M|$ vertices of
degree $1$ and $2$ is obviously finite (here $|M|$ denotes the
number of elements in $M$). So there is only a finite number of
possible types for a filling of the space~$\cal M$ and  for each
type there exists a generalized  minimal parametric filling by the
previous lemma. Now we obtain a generalized minimal filling of the
space~$\cal M$ by choosing the filling of minimal weight from this
finite collection of generalized minimal parametric fillings.

Suppose the filling has some vertices of degree different from $1$
and $3$. Then we split each of the vertices by adding  edges of zero
weight. Thus we obtain a binary tree. This operation does not change
the graph weight and the path lengths, so  the resulting generalized
filling is still minimal. This completes the proof.

Let  $\cal G$ be a generalized filling of the space $\cal M$. A path
$\gamma$ in $\cal G$ connecting points $x$ и $y$ is called {\it
exact} if  $x, y \in M$ and its weight is equal to the distance
between its endpoints in the space $\cal M$, i.e.
$w(\gamma)=\rho(x,y)$. We denote by ${\rm deg}(v)$ the degree of the
vertex $v$, i.e. the number of edges incident to it.

\begin{theorem}\label{lem: exact paths}
\ \\ {\noindent\bf 1.} For any edge in a generalized minimal
parametric filling there exists an exact path containing the edge.

{\noindent\bf 2.} For any two adjacent edges in a generalized
minimal  filling there exists an exact path containing the edge.

{\noindent\bf 3.} For any pair of edges incident to an interior
vertex of degree~$3$ in a generalized minimal parametric filling
there exists an exact path containing the pair of edges.

{\noindent \bf 4.} Fix an interior vertex~$v$ in a generalized
minimal parametric filling and a subset of $m$ edges incident to the
vertex. Suppose \mbox{$m>\frac{1}{2}{\rm deg}(v)$}. Then there
exists an exact path containing a pair of edges from the subset.

\end{theorem}
{\bf Proof.} {\bf 1.} The first part is clear since otherwise we
could decrease the edge weight and obtain a filling with smaller
weight.

{\bf 2. } Suppose  there are no exact paths containing the two
adjacent edges $xv$ and $vy$ in a generalized filling ${\cal G}=
(V,E,w)$. Consider the graph $G'$ with one additional vertex: $V' =
V\cup \{u\}$ and with the edge set $E' =
E\setminus\{vx,vy\}\cup\{ux,uy,uv\}$. Define $w'(uv)=\varepsilon$,
$w'(xu)=w(xv)-\varepsilon$, $w'(yu)=w(yv)-\varepsilon$, $w'(e)=w(e)$
for other edges. Then the lengths of all the paths containing the
two edges $xv$ and $vy$ are reduced  by $2\varepsilon$ and the
lengths of other paths with boundary endpoints  remain unchanged.
Since no exact paths contained the two edges $xv$ and $vy$,
$\varepsilon>0$ can be chosen in such a way  that $(G',w')$ is a
filling. But $w'({\cal G}')=w({\cal G})-\varepsilon$, which
contradicts minimality of the original filling~$\cal G$.

{\bf 3.} Follows from  {\bf 4} by putting $m=2, {\rm deg}(v)=3$.

{\bf 4.} Suppose  there is a set of $m>\frac{1}{2}{\rm deg}(v)$
edges incident to an interior vertex~$v$ provided that no two edges
from the set are contained in an  exact path. Define the weight
function $w'$ on $G$ by reducing the weight of every edge from the
set by $\varepsilon$ and by adding $\varepsilon$ to the weight of
every other edge incident to~$v$, while the remaining edge weights
are left unchanged. By the assumption $m>{\rm deg}(v)-m$, so for any
$\varepsilon>0 $ the weight of the graph decreases. On the other
hand, for any  sufficiently small  $\varepsilon>0$ the graph
$(G,w')$ is a filling. This contradiction completes the proof.

\medskip

\noindent{\bf Remark 5.} Each of the last three lemmas is a
straightforward generalization of the corresponding result in the
case of non-generalized fillings, see~\cite{it}. Moreover, the proof
in the case of generalized fillings appears to be simpler since one
does not have to check if the edge weights are non-negative.

\begin{theorem}\label{lem: boundary edge weight is positive}
Consider a boundary vertex in a generalized minimal parametric
filling and consider an edge incident to the vertex. Then the edge
weight is non-negative.
\end{theorem}
{\bf Proof.} Let $p$ be a boundary  vertex and $px$ be an edge
incident to it. Suppose $x$ is boundary. Then
\mbox{$w(px)\geqslant\rho(p,x)\geqslant 0$}, and the proof is
finished. In the converse case, the degree of vertex $x$ is not less
than~$3$ and by lemma~\ref{lem: exact paths}, there exists an exact
path containing a pair of edges incident to $x$ and different from
the edge $px$. Denote the endpoints of the path by $q$ and $r$ (then
we denote the path by $q$-$r$). Consider the paths $p$-$q$ and
$p$-$r$. These two paths have only one common edge $px$. The path
$q$-$r$ consists of all the edges contained in paths $p$-$r$ and
$p$-$q$ excluding the edge $px$. Therefore:
$$ \rho(q,r)= d_w(q,r) = d_w(q,p)+d_w(r,p)-2w(px) \geqslant
\rho(q,p)+\rho(r,p)-2w(px)\mbox{,}$$

\noindent on the other hand,
$\rho(q,p)+\rho(r,p)\geqslant\rho(q,r)$, so
 $w(px)\geqslant0$. This completes the proof.

\medskip

\noindent{\bf Remark 6.} Suppose the finite space under
consideration is a non-degenerate metric space, i.e. the strict
inequalities $\rho(x,y)>0$ and $\rho(x,y)+\rho(y,z)>\rho(x,z)$ hold
for any $x\ne y\ne z\ne x$. Then the weight of any edge incident to
a boundary vertex in a generalized minimal parametric filling of
this space is strictly positive. This easily follows from the proof
of lemma~\ref{lem: exact paths}.

\begin{wrapfigure}[14]{r}{100pt}
\includegraphics[scale = 0.7]{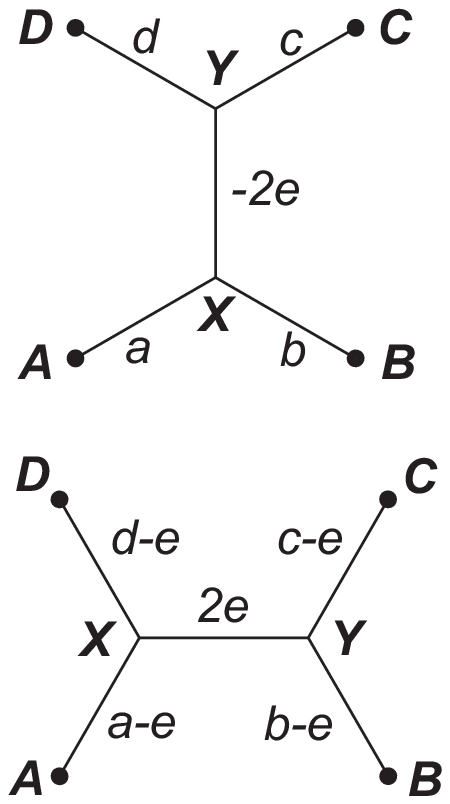}
\end{wrapfigure}

\medskip

\noindent \textbf{Modification.} Let ${\cal G } =(V,E,w)$ be a
generalized filling of space $\cal M$, and let $X$ and $Y$ be
vertices of degree $3$ connected by a negative weight edge:
$w(XY)=-2e<0$. Denote by $A$ and $B$ the vertices adjacent to~$X$
and different from~$Y$ and denote by $C$ and $D$ the vertices
adjacent to  $Y$ and different from~$X$. Denote $w(XA)=a$,
$w(XB)=b$, $w(YC)=c$, $w(YD)=d$. We construct the modified tree
$G'=(V',E')$ as follows.

The vertex set is the same: $V'=V$. The edge set is modified: $E' =
E\setminus \{XB, YD \}\cup\{XD,YB\}$. Define the weight function:
$w'(XY)=2e$, $w'(XA)=a - e$, $w'(YB)=b - e$, $w'(YC)=c - e$,
$w'(XD)=d - e$, $w'=w$ on other edges of $G'$.

\begin{theorem}\label{lem: modification}
\ \\
\noindent{\rm\textbf{1.}} The modified graph ${\cal G}'=
(G',w')$ is a generalized filling of $\cal M$ with the same weight
$w'({\cal G'})=w(G)$.

\noindent{\rm \textbf{2.}} Suppose $u,v\in {M}$. Then
$d_{w'}(u,v)\geqslant d_w(u,v)$.

\noindent{\rm \textbf{3.}} Suppose a path in $G$ with endpoints
$u,v\in {M}$ contained the edge sequence $AX$, $XY$, $YC$. Then the
path length has increased in $G'$, i.e. $d_{w'}(u,v)
> d_w(u,v)$.
\end{theorem}

\noindent{\bf Proof.} Statement {\bf 1} follows from  {\bf 2}. The
fact $w'({\cal G'})=w({\cal G})$ is obvious. Let us prove
statements {\bf 2} и {\bf 3}. Consider a path with boundary
endpoints. Suppose the path contains neither of the vertices $X$
and $Y$. Then the path length obviously does not change after
modification since the whole path remains unchanged. Suppose the
path in $G$ contained the edge sequence $AX$, $XB$ (with weight
$a+b$). Then the corresponding path in $G'$ contains the edge
sequence $AX$, $XY$, $YB$ (with the same weight $a-e+2e+b-e =
a+b$). Therefore in this case the path length also remains
unchanged after modification. The situation is similar in the
cases where the path in $G$ contained the edge sequence $CY$,
$YD$, or $BX$, $XY$, $YC$, or $AX$, $XY$, $YD$.

 Now suppose the path in $G$ contained the edge sequence
$AX$, $XY$, $YC$. Then the path length increases after modification
since an arc of weight $a-2e+d$ is replaced by an arc of weight
$a-e+2e+d-e=a+d$, which is strictly greater than $a-2e+d$. The case
of $BX$, $XY$, $YD$ is similar. The proof is finished.

\Name{Proof of the main theorem}\label{sec: the-main-theorem}

Let us recall the main theorem.

\medskip

\noindent {\bf Theorem.} {\it Let ${\cal M}$ be a finite
pseudometric space. Then ${\rm mf_{-}}({\cal M})={\rm mf}({\cal
M})$.}

\medskip

\noindent{\bf Proof.} As we have already mentioned above, it is
obvious that ${\rm mf}_{-}({\cal M})\leqslant{\rm mf}({\cal M})$.
Hence it is sufficient to prove that there exists a filling with
weight  ${\rm mf}_{-}({\cal M})$ and without negative edge
weights.

Choose a filling $\cal G$ with the minimum number of exact paths
among all the generalized minimal fillings of $\cal M$ that are
binary trees (at least one such filling exists by lemma~\ref{lem:
mf- exists}). Let us prove that there are no negative edge weights
in filling $\cal G$ and thus $G$ provides an example of
non-generalized filling with weight~${\rm mf}_{-}({\cal M})$.

Suppose there is an edge of negative weight. Denote it by $XY$. By
lemma~\ref{lem: exact paths}, there exists an exact path~$\gamma$
containing the edge $XY$. The vertices $X$ and $Y$ are interior
since otherwise the weight of $XY$ would be non-negative by
lemma~\ref{lem: boundary edge weight is positive}. So we can modify
the graph as described above. Moreover, by lemma~\ref{lem:
modification} we can modify the graph in such a way that the weight
of the path $\gamma$ increases and the weights of other paths with
boundary endpoints do not decrease. Since the modification preserves
the graph weight the modified graph is a generalized minimal filling
of $\cal M$. But the number of exact paths has reduced by at least 1
in the modified graph. This contradicts the choice of~$\cal G$. The
proof is finished.

\medskip

\noindent{\bf Remark 7.} The theorem does not hold if the function
$\rho$ on the boundary set~$M$ violates the triangle inequality. To
be exact, only one of our lemmas uses the triangle inequality. This
is lemma~\ref{lem: boundary edge weight is positive} on
non-negativity of boundary edge weights.

\medskip

\noindent{\bf Example 3.} Consider the set $\{x,y,z\}$ and the
function \mbox{$\rho\colon \rho(x,y)=1$}, $\rho(y,z)=2$,
$\rho(z,x)=5$. The triangle inequality does not hold but still one
can search for the minimal filling weight. It is easy to see that
in this case the generalized minimal filling weight and the
minimal filling weight are different. They are $4$ and $5$
respectively. And there is  a boundary edge of negative weight in
the generalized minimal filling.

\end{document}